\documentclass[leqnographics,dvips,emlines, openbib,
    keyval,epsfig,makeidx,amsfonts,amssymb]{article}
\usepackage{amssymb}
\usepackage{amsfonts}
\usepackage{amsmath}
\usepackage{amsxtra}
\title{A Note on Doubly Warped Product Contact $CR$-Submanifolds
in trans-Sasakian Manifolds}
\author{ Marian-Ioan Munteanu
\thanks{ Beneficiary of a Research Fellowship of the Ministerie van de
Vlaamse Gemeenschap,\newline ref. 13DA-IZ/JDW/ROE/2004-2005/04-2439
      at the Katholieke Universiteit Leuven, Belgium}
}
\def\proof{{\sc Proof. }}

\newtheorem{contor}{1.}

\newtheorem{theorem}[contor]{Theorem}

\newtheorem{corollary}[contor]{Corollary}

\def\proof{{\sc Proof.\ }}
\newcommand{\gata}{\hfill\hskip -1cm \rule{.5em}{.5em}}

\date{ }

\begin{document}

\maketitle

\begin{abstract}

\noindent
Warped product $CR$-submanifolds in K\"ahlerian manifolds were intensively
studied only since 2001 after the impulse given by B.Y. Chen in \cite{kn:BYC133},
\cite{kn:BYC134}. Immediately after, another line of research, similar to that
concerning Sasakian geometry as the odd dimensional version of K\"ahlerian
geometry, was developed, namely warped product contact $CR$-submanifolds in
Sasakian manifolds (cf. \cite{kn:IM}, \cite{kn:MIM}). In this note we proved
that there exists no proper doubly warped product contact $CR$-submanifolds
in trans-Sasakian manifolds.

\vspace{2mm}


\vspace{2mm}


\end{abstract}

\section{ Preliminaries}

Let $M$ be a Riemannian manifold isometrically immersed in a trans-Sasakian
manifold $(\widetilde M,\phi,\xi,\eta,\widetilde g)$. M is called a {\em contact
$CR$-submanifold} if there exists on $M$ a differentiable $\phi$-invariant
distribution ${\mathcal{D}}$ (i.e. $\phi_x{\mathcal{D}}_x\subset {\mathcal{D}}_x$,
for all $x\in M$) whose orthogonal complement ${\mathcal{D}}^\bot$ of
${\mathcal{D}}$ in $T(M)$ is a $\phi$-anti-invariant distribution on $M$,
i.e. $\phi_x{\mathcal{D}}^\bot_x\subset T(M)_x^\bot$ for all $x\in M$. Here
$T(M)^\bot\longrightarrow M$ is the normal bundle of $M$.

\vspace{2mm}

It is customary to require that $\xi$ be tangent to $M$ than normal
(in case when $\widetilde M$ is Sasakian, by Prop. 1.1,
K.Yano \& M.Kon \cite{kn:KYMK}, p.43, in this later case,
$M$ must be $\phi$-anti-invariant, i.e. $\phi_xT_x(M)\subset T(M)^\bot_x$,
$\forall x\in M$).

\vspace{2mm}

Given a contact $CR$-submanifold $M$ of a trans-Sasakian manifold $\widetilde M$,
either $\xi\in{\mathcal{D}}$, or $\xi \in {\mathcal{D}}^\bot$
(the result holds only due the almost contact structure conditions).
Therefore, the tangent space at each point decomposes orthogonally as
\begin{equation}
T(M)=H(M)\oplus{\mathbf{R}}\xi\oplus E(M)
\end{equation}
where $\phi H(M)=H(M)$, $\phi^2{}_{|_{H(M)}}=-I_{H(M)}$
($H(M)$ is called {\em the Levi distribution} of $M$)
and $\phi E(M)\subset T(M)^\bot$. We will consider only {\em proper}
contact $CR$-submanifolds ($\dim H(M)>0$), i.e. $M$ is neither
$\phi$-invariant, nor $\phi$-anti invariant. Remark that both
${\mathcal{D}}=H(M)$ and ${\mathcal{D}}=H(M)\oplus{\mathbf{R}}\xi$
organize $M$ as a contact $CR$-submanifold, but $H(M)$ is never
integrable if $d\eta$ is non-degenerate (for example if the ambient is
a contact manifold, e.g. a Sasakian manifold).

\vspace{2mm}

An {\em almost contact structure} $(\phi, \xi, \eta)$
on a $(2m+1)$-dimensional manifold $\widetilde M$,
is defined by $\phi\in {\mathcal{T}}^1_1(\widetilde M)$,
$\xi\in\chi(\widetilde M)$, $\eta\in\Lambda^1(\widetilde M)$ satisfying
the following properties: $\phi^2=-I+\eta\otimes\xi$, $\phi\xi=0$,
$\eta\circ\phi=0$, $\eta(\xi)=1$.
If on $\widetilde M$ we have the metric $\widetilde g$, then the compatibility
condition is required, namely,
$\widetilde g(\phi X, \phi Y)=\widetilde g(X,Y)-\eta(X)\eta(Y)$.

\vspace{2mm}

An almost contact metric structure $(\phi, \xi, \eta, \widetilde g)$
on $\widetilde M$ is called a {\em trans-Sasakian structure} if
$(\widetilde M\times {\mathbf{R}}, J, G)$ belongs to the class ${\mathcal{W}}_4$
of the Gray-Hervella classification of almost Hermitian manifolds
(see \cite{kn:GH}), where $J$ is the almost complex structure on the product
manifold $\widetilde M\times {\mathbf{R}}$ defined by
$$
J\left(X, f\frac d{dt}\right)=\left(\phi X-f\xi,\eta(X)\frac d{dt}\right)
$$
for all vector fields $X$ on $\widetilde M$ and smooth functions $f$ on
$\widetilde M\times {\mathbf{R}}$. Here $G$ is the product metric on
$\widetilde M\times {\mathbf{R}}$. If $\widetilde\nabla$ denotes
the Levi Civita connection (of $\widetilde g$), then the following
condition
\begin{equation}
\label{nablaphi}
(\widetilde\nabla_X\phi)Y=\alpha (g(X,Y)\xi-\eta(Y)X)+
    \beta (g(\phi X,Y)\xi-\eta(Y)\phi X)
\end{equation}
holds for some smooth functions $\alpha$ and $\beta$ on $\widetilde M$.
We say that the trans-Sasakian structure is {\em of type} $(\alpha,\beta)$.
From the formula (\ref{nablaphi}) it follows that
\begin{equation}
\widetilde\nabla_X\xi=-\alpha \phi X+\beta(X-\eta(X)\xi).
\end{equation}
(See e.g. \cite{kn:DEB01} for details.)

\vspace{2mm}

We note that trans-Sasakian structures of type $(0,0)$ are
cosymplectic, trans-Sasakian structures of type $(0, \beta)$
are $\beta$-Kenmotsu and a trans-Sasakian structure of type
$(\alpha, 0)$ are $\alpha$-Sasakian.

\vspace{2mm}

For any $X\in\chi(M)$ we put $PX=tan(\phi X)$ and $FX=nor(\phi X)$,
where $tan_x$ and $nor_x$ are the natural projections associated to
the direct sum decomposition
$T_x(\widetilde M)=T_x(M)\oplus T(M)_x^\bot,$ $x\in M$.
We recall the Gauss formula
\begin{equation}
\widetilde\nabla_XY=\nabla_XY+h(X,Y)
\end{equation}
for any $X,Y\in\chi(M)$. Here $\nabla$ is the induced connection and
$h$ is the second fundamental form of the given immersion. Since $\xi$
is tangent to $M$ we have
\begin{equation}
\label{eq5}
P\xi=0,\quad F\xi=0,\quad \nabla_X\xi=-\alpha PX+\beta (X-\eta(X)\xi),
\quad h(X,\xi)=-\alpha FX,\ X\in\chi(M).
\end{equation}

\vspace{2mm}

\section{Doubly Warped Products}

Doubly warped products can be considered as generalization of warped
products. A {\em doubly warped product} $(M,g)$ is a product manifold
of the form $M=_fB\times_bF$ with the metric $g=f^2g_B\oplus b^2g_F$,
where $b:B\longrightarrow (0,\infty)$ and
$f:F\longrightarrow (0,\infty)$ are smooth maps and $g_B, g_F$
are the metrics on the Riemannian manifolds $B$ and $F$ respectively.
(See for example \cite{kn:BU01}.) If either $b\equiv1$
or $f\equiv1$, but not both, then we obtain a ({\em single})
{\em warped product}. If both $b\equiv1$ and $f\equiv1$, then we have a
{\em product manifold}. If neither $b$ nor $f$ is constant, then we have a
{\em non trivial doubly warped product}.

\vspace{2mm}

If $X\in\chi(B)$ and $Z\in\chi(F)$, then the Levi Civita connection is
\begin{equation}
\label{patru}
\nabla_XZ=Z(\ln f)X+X(\ln b)Z.
\end{equation}

\vspace{2mm}

We give the main result of this section. Namely we have

\vspace{2mm}

\begin{theorem}
There is no proper doubly warped product contact $CR$-submanifolds
in trans-Sasakian manifolds.
\end{theorem}
\proof
Let $M=_{f_2}N^\top\times_{f_1}N^\bot$ be a doubly warped product
contact $CR$-submanifold in a trans-Sasakian manifold
$(\widetilde M, \phi, \xi, \eta, \widetilde g)$,
i.e. $N^\top$ is a $\phi$-invariant submanifold and
$N^\bot$ is a $\phi$-anti-invariant submanifold.

\vspace{.5mm}

As we have already seen $\xi\in{\mathcal{D}}$ or $\xi\in{\mathcal{D}}^\bot$.

\vspace{.5mm}

Case 1. $\xi\in{\mathcal{D}}$ i.e. $\xi$ is tangent to $N^\top$. Taking
$Z\in\chi(N^\bot)$ we have $\nabla_Z\xi=\beta Z$. On the other hand, from
(\ref{patru}), we get
$$
\nabla_\xi Z=Z(\ln f_2)\xi+\xi(\ln f_1) Z.
$$
It follows, since the two distributions are orthogonal, that
\begin{equation}
\left\{\begin{array}{l}
\xi(\ln f_1)=\beta,\\
Z(\ln f_2)=0\quad {\rm for\ all\ } Z\in{\mathcal{D}}^\bot.
\end{array}\right.
\end{equation}
The second condition yield $f_2\equiv constant$, so, we cannot have
doubly warped product contact $CR$-submanifolds of the form
${}_{f_2} N^\top\times_{f_1} N^\bot$, with $\xi$ tangent to $N^\top$,
other than warped product contact $CR$-submanifolds. Moreover, in this
case $\beta$ is a smooth function on $N^\top$.

\vspace{1mm}

Case 2. $\xi\in{\mathcal{D}}^\bot$, i.e. $\xi$ is tangent to $N^\bot$.
Taking $X\in\chi(N^\top)$ we have
$\nabla_X\xi=-\alpha PX+\beta X$ in one hand and
$\nabla_X\xi=X(\ln f_1)\xi+\xi(\ln f_2)X$ in the other hand. Since the
two distributions are orthogonal, we immediately get
\begin{equation}
\label{sase}
\left\{\begin{array}{l}
\xi(\ln f_2)X=-\alpha PX+\beta X,\\
X(\ln f_1)=0\quad {\rm for\ all\ } X\in{\mathcal{D}}.
\end{array}\right.
\end{equation}
The second condition in (\ref{sase}) shows that $M$ is a $CR$-warped product
between a $\phi$-anti invariant
manifold $N^\bot$ tangent to the structure vector field $\xi$
and an invariant manifold $N^\top$.
If $\dim{\mathcal{D}}=0$, then $M$ is a $\phi$-anti invariant submanifold
in $\widetilde M$. Otherwise, one can choose $X\neq0$, and thus
$X$ and $PX$ are linearly independent.
Using first relation in (\ref{sase}) one gets $\alpha=0$ and $\beta=\xi(\ln f_2)$.
This means that the ambient manifold is $\beta$-Kenmotsu with
$\beta\in C^\infty(N^\bot)$.

This ends the proof.
\gata

\vspace{2mm}

In 1992, J.C. Marero in (\cite{kn:JCM}) showed that a trans-Sasakian manifold
of dimension $\geq 5$ is either $\alpha$-Sasakian, $\beta$-Kenmotsu or
cosymplectic. So, we can state

\begin{corollary} \rm
Let $\widetilde M$ be

\begin{enumerate}

\item{ }  {\em either an $\alpha$-Sasakian manifold},

\item{ } {\em or a $\beta$-Kenmotsu manifold},

\item{ } {\em or a cosymplectic manifold}.

\end{enumerate}

Then, there is no proper doubly
warped product contact $CR$-submanifolds in $\widetilde M$.
More precisely we have,

$\checkmark\ $ if $\xi\in{\mathcal{D}}$:

\qquad \ \ $M=\widetilde N^\top\times {}_fN^\bot$,
$\xi$ is tangent to $N^\top$ and $f\in C^\infty(N^\top)$, where
$\widetilde N^\top$ is the manifold $N^\top$ with a homothetic metric
$c^2g_{N^\top}$, $(c\in{\mathbf{R}})$.
Moreover, in case 2, $\beta$ is a smooth funtion on $N^\top$.

$\checkmark\ $ if $\xi\in{\mathcal{D}}^\bot$:

\qquad \ \ 1. $M$ is a $\phi$-anti-invariant submanifold in $\widetilde M$
                ($\dim {\mathcal{D}}=0$);

\qquad 2-3. $M=\widetilde N^\bot\times {}_fN^\top$,
$\xi$ is tangent to $N^\bot$ and $f\in C^\infty(N^\bot)$, where
$\widetilde N^\bot$ is the manifold $N^\bot$ with a homothetic metric
$c^2g_{N^\bot}$, $(c\in{\mathbf{R}})$. Moreover, in case 2,
$\beta$ is a smooth function on $N^\bot$
\end{corollary}
\proof
The statements follow from the Theorem 1.
\gata

\vspace{2mm}

It follows that there is no warped product contact $CR$ submanifolds
in $\alpha$-Sasakian and cosymplectic manifolds
in the form $N^\bot\times {}_fN^\top$ (with $\xi$ tangent to $N^\top$)
other than product manifolds ($f$ must be a real positive constant)
-- after a homothetic transformation of the metric on $N^\top$.
If the ambient is $\beta$-Kenmotsu, then

\qquad -- there is no contact $CR$-product submanifolds

\qquad -- there is no contact $CR$-warped product submanifolds of type
$N^\bot\times{}_fN^\top$ (with $\xi$ tangent to $N^\top$)
(see also \cite{kn:MAEM}).

\vspace{3mm}

We will give an example of warped product contact $CR$-submanifold
of type $N^\bot\times {}_fN^\top$ in Kenmotsu manifold, with $\xi$
tangent to $N^\bot$:

\vspace{2mm}

Consider the complex space ${\mathbf{C}}^m$ with the usual Kaehler
structure and real global coordinates $(x^1,y^1,\ldots,x^m,y^m)$.
Let $\widetilde M={\mathbf{R}}\times {}_f{\mathbf{C}}^m$ be the warped
product between the real line ${\mathbf{R}}$ and ${\mathbf{C}}^m$,
where the warping function is $f=e^z$, $z$ being the global
coordinate on $\mathbf{R}$. Then, $\widetilde M$ is a Kenmotsu manifold
(see e.g. \cite{kn:DEB01}). Consider the distribution
$\mathcal{D}={\rm span\ }\left\{
\frac\partial{\partial x^1},\frac\partial{\partial y^1},
\ldots, \frac\partial{\partial x^s},
\frac\partial{\partial y^s}\right\}$
and
$\mathcal{D}^\bot={\rm span\ }\left\{
\frac\partial{\partial z},
\frac\partial{\partial x^{s+1}},
\ldots,\frac\partial{\partial x^m}\right\}$
which are integrable and denote by $N^\top$ and $N^\bot$
the integral submanifolds, respectively. Let
$g_{N^\top}=\sum\limits_{i=1}^s\left((dx^i)^2+(dy^i)^2\right)$
and $g_{N^\bot}=dz^2+e^{2z}\sum\limits_{a=s+1}^m(dx^a)^2$
be Riemannian metrics on $N^\top$ and $N^\bot$, respectively.
Then, $M=N^\bot\times{}_fN^\top$ is a contact $CR$-submanifold,
isometrically immersed in $\widetilde M$. Here, the warping function
is $f=e^z$.

\vspace{2mm}

{\bf Acknowledgements.} This work was done while the author was visiting
the Department of Mathematics, Catholic University of Leuven, Belgium,
in the period February $-$ May 2005. I wish to thank Prof. Franki Dillen
for his helpful advices and for his warm hospitality during my stay in Leuven.

\vspace{2mm}

\scriptsize { Permanent address:
\hspace{60mm}  Temporary address:

\sc

University 'Al.I.Cuza' of Ia\c si,
\hspace{47mm}   Katholieke Universiteit Leuven,

Faculty of Mathematics,
\hspace{53mm}    Departement Wiskunde

Bd. Carol I, no.11,
\hspace{60mm}   Celestijnenlaan 200B

700506 Ia\c si,
\hspace{70mm}       B-3001 Leuven (Heverlee)

Romania
\hspace{74mm}    BELGIUM }

\vspace{4mm}

e-mail: munteanu@uaic.ro,               \hspace{70mm}  Received,

$\qquad\quad$ munteanu2001@hotmail.com    \hspace{55mm} January 20, 2006

\normalsize

\end{document}